\documentclass[11pt,roman]{article}

\usepackage{color} % To make text of different colors
\usepackage{graphicx}% Include figure files

\begin{document}

% DO NOT INSERT PAGE NUMBERS
\pagestyle{empty}

% DO NOT PRINT DATE
\date{}

% MAKE TITLE BOLD 14 PT
\vspace*{13pt}\begin{center} \Large\textbf{An Improved Procedure for
Selecting the Profiles of Perfectly Matched Layers}
\end{center}

% PUT AUTHORS' NAMES HERE
\vspace*{13pt}\begin{center} \large\textbf{Jiawei~Zhang$^1$}
\end{center}

% PUT AFFILIATION HERE
\begin{center} % affiliation 1
$^1$ Department of Mathematics, Zhejiang~University, Hangzhou, 310058, China. \\
\textcolor{blue} {\underline {jiaweiz@zju.edu.cn}}
\end{center}

\setlength{\lineskip}{12pt}
% TYPE YOUR ABSTRACT HERE
\vspace*{26pt}\noindent\textbf{Abstract:} The perfectly matched
layers (PMLs), as a boundary termination over an unbounded spatial
domain, are widely used in numerical simulations of wave propagation
problems. Given a set of discretization parameters, a procedure to
select the PML profiles based on minimizing the discrete
reflectivity is established for frequency domain simulations. We, by
extending the function class and adopting a direct search method,
improve the former procedure for traveling waves.\\

% TYPE YOUR KEYWORDS HERE
 \vspace*{13pt}

% define Section counter
\newcounter{Section}
% define Subsection counter
\newcounter{Subsection}

% START A NEW SECTION HERE
% increment counter
\addtocounter{Section}{1}
% ENTER SECTION TITLE AFTER DOT
\begin{center} \textbf{\arabic{Section}. Introduction}
\end{center}

% TYPE THE TEXT OF YOUR 1ST SECTION HERE
Perfectly matched layers (PMLs)~\cite{Berenger1} are widely used as
boundary terminations in problems to be solved over an unbounded
spatial domain. Physically, this is an approach to surround a
numerical problem domain with a layer of a material which creates as
little numerical reflection as possible, while also attenuating
waves that enter from the problem interior. Mathematically, when a
PML is used to truncate the $x$ axis, $x$ is actually replaced by
$\hat{x}=\int^{x}{(1+i\sigma(\tau))d\tau}$, where $\sigma$ is a real
function satisfying certain conditions.

Nevertheless, in actual numerical simulations, due to the finite
thickness of the PML, reflections occur when plane waves incidence
upon the PML. Note that PML is generally problem-dependent, for
example, the reflection is dependent on the discretization scheme.
For frequency domain simulations, based on minimizing the average
discrete reflectivity, Ya Yan Lu~\cite{Lu1} gave a practical
procedure for selecting the optimal PML profile $\sigma$ where the
function class is restricted to simple powers. In this paper, to
further reduce the average reflectivity for traveling waves, we
consider a rational function class for $\sigma$; to avoid
time-consuming computations, we use the Nelder-Mead simplex method
to determine the coefficients of $\sigma$. Besides, we also give a
simpler $\sigma$ which is easy to optimize. Numerical simulations
demonstrate that our improved procedure is much better and
practical.

% START A NEW SECTION HERE
% increment counter
\addtocounter{Section}{1}
% ENTER SECTION TITLE AFTER DOT
\vspace*{13pt}
\begin{center} \textbf{\arabic{Section}. Motivations to Improve the Former Procedure}
\end{center}

% TYPE THE TEXT OF YOUR 2TH SECTION HERE

For traveling waves, we consider a two-dimensional waveguiding
structure in ~\cite{Lu1}. The $y-$component of the electric field of
a transverse wave satisfies a Helmholtz equation. Consider that the
structure is unbounded in the negative $x$ direction and the medium
is homogeneous (refractive index $n \equiv n_{0}$) for $x < G$, we
then truncate the negative $x$ axis by a PML. For $D < H < G$, we
define $\sigma(x)$ such that $\sigma(x)=0$ for $x \geq H$ and
$\sigma(x)>0$ for $x < H$ while $\sigma(H)=0$ and $\sigma'(H)=0$.
Replacing $x$ by $\hat{x}=\int^{x}{(1+i\sigma(t))dt}$ gives rise to
\begin{equation}
    s^{-1}\partial_{x}({s^{-1}\partial_{x}{u}})+
    \partial^{2}_{z}{u}+k_{0}^{2}n_{0}^{2}u=0
\end{equation}
where $k_{0}$ is the free space wavenumber, and the time dependence
is $e^{-i\omega t}$. At $x=D$, we use a simple zero boundary
condition: $u=0$. The actual PML is the layer $D<x<H$. For $H<x<G$,
(1) has a plane wave solution
\begin{equation}
    u=e^{i(-\alpha x + \beta z)} + R e^{i(\alpha x + \beta z)}
\end{equation}
where the second term is the reflected wave due to the incidence of
plane waves upon the PML with reflection coefficient $R$. When $x$
is discretized, $R$ depends on $\sigma$~\cite{YYS1}. By a
second-order finite difference approximation in the transverse
direction $x$~\cite{Lu1}, we can easily find $R$ exactly by solving
a linear equation system.

Thus, to select the optimal PML profile $\sigma$ is to find a
$\sigma$ such that the following $\overline{|R|}$ is minimized:
\begin{equation}
    \overline{|R|} = \frac{2}{\pi}\int_0^{\pi/2}|R(\theta)|d\theta.
\end{equation}
For convenience, let
\begin{equation}
    \tau(x)=\frac{x-H}{D-H}.
\end{equation}
In Lu's procedure~\cite{Lu1}, since $\sigma$ is limited to be a
simple power, i.e.
\begin{equation}
    \sigma = S\tau^p,
\end{equation}
we only need to determine $p$ and the dimensionless scaling
parameter $S$ such that $\overline{|R|}$ is minimized. The optimal
values of $S$ were computed for $p=2,...,5$ respectively, among
which the one giving the least $\overline{|R|}$ was chosen to be the
overall optimal PML profile.

In fact, the numerical result can be more satisfying if we give up
the restriction of $\sigma$ to be simple powers. In practice, people
usually use $\sigma = \frac{S\tau^3}{1+\tau^2}$ for the PML
profile~\cite{LZ1}. This gives us a start to find a better profile.

After numerical computation and a little adjustment, $\sigma =
\frac{S\tau^3}{1+\tau}$ turns out to be better in this situation. So
in this paper, we first investigate into this rational function
class:
\begin{equation}
    \sigma = \frac{a_2\tau^2+a_3\tau^3+...+a_p\tau^p}
    {1+\tau}, a_p>0, p\geq2
\end{equation}
where its coefficients are to be optimized to give a minimal value
of the average discrete reflectivity $\overline{|R|}$. Due to the
condition $\sigma>0$, we let $a_p>0$ for simplicity while losing
certain generality.

However, if the PML profile is defined by such rational function,
the work to determine the optimal values of its coefficients becomes
much more time-consuming. To save time, we adopt the Nelder-Mead
(NM) simplex method due to the fact that the problem is nonlinear
and the derivative information of $\overline{|R|}$ is unavailable.

It is natural to ask why we choose this method. We offer two
answers. First, in the NM method, the objective function is
evaluated at the vertices of a simplex, and movement is away from
the poorest value. This method tends to work so well in practice by
producing a rapid initial decrease in function values. Second, when
we consider a simpler $\sigma$ with two parameters, the NM algorithm
gives answer in a much shorter time. Thus, due to its powerful local
descent property, we decide to adopt this method though it may not
give globally optimized solution.~\cite{LRWW1}

Note that the NM method is for unconstrained problems, but $a_k>0$
in (6) is actually a constrainment. To overcome this, we define
$\sigma$ as below
\begin{equation}
    \sigma = \frac{|a_2|\tau^2+|a_3|\tau^3+...+|a_p|\tau^p}
    {1+\tau}, a_p>0, p\geq2.
\end{equation}
When similar situation occurs, we will tackle it in this way again
without further remarks.

Furthermore, after optimization for (6), we try to define $\sigma$
by a simpler function class which largely conserves the good
properties of the former one but is much easier to optimize and more
useful in practice.

% START A NEW SECTION HERE
% increment counter
\addtocounter{Section}{1}
% ENTER SECTION TITLE AFTER DOT
\vspace*{13pt}\begin{center} \textbf{\arabic{Section}. The Improved
Procedure with Numerical Results}
\end{center}

% START 1ST SUBSECTION HERE
% increment counter
\addtocounter{Subsection}{1}
% ENTER SUBSECTION TITLE AFTER DOT
\begin{flushleft} \textit{\Alph{Subsection}. General Results}
\end{flushleft}\vspace*{-7pt}

% TYPE THE TEXT OF YOUR 1ST SUBSECTION HERE
In numerical simulations, we adopt the conditions in \cite{Lu1} so
that we could compare the average discrete reflectivity
$\overline{|R|}$ between the two procedures. Explicitly we have the
wavelength $\lambda_0=1$ $\mu$ m, $k_0=2\pi /\lambda_0$, $n_0=1$.

Consider an example in~\cite{Lu1}: A PML with thickness of five
grids ($m=5$), where the grid size $h=\frac{1}{20}\lambda_0=0.05$
$\mu$ m. Under the restriction on $\sigma$ in (5): $\sigma =
S\tau^p, p\le 5$, the author obtained an optimal profile: $p=3$ and
$S=100.4$, which gave rise to $\overline{|R|}=0.013$.

Through our improved procedure, we obtain a simple and better
result: $\sigma = (23.6\tau^2+35.9\tau^5)/(1-\tau)$ which gives rise
to $\overline{|R|}$ = 0.0047, only 36\% of the former one. Moreover,
if we give up the simplicity and allow a higher order of the
rational function, we could further get an result of 0.0031, only
24\% of the former one.

\begin{table}[!h]
\tabcolsep 0pt \caption{Local optimal values of the coefficients of
$\sigma$ defined by (6). Results are derived by the Nelder-Mead
simplex method. At $p=6, 9, 12$, the algorithm stops because the
number of function evaluation has exceeded the preset limit of 2000.
The actual number of iteration will be bigger if we increase the
limit.} \vspace*{-12pt}
    \begin{center}
    \def\temptablewidth{0.9\textwidth}
    {\rule{\temptablewidth}{1pt}}
        \begin{tabular*}{\temptablewidth}{@{\extracolsep{\fill}}cccccccccccccc}
            $p$ &$a_2$ &$a_3$ &$a_4$ &$a_5$ &$a_6$ &$a_7$ &$a_8$ &$a_9$ &$a_{10}$ &$a_{11}$ &$a_{12}$ &Iterations &$\overline{|R|}$ \\   \hline
            2 &74.2  &  &  &  &  &  &  &  &  &  &  & 21 &  0.019  \\
            3 &38.2  &108.7  &  &  &  &  &  &  &  &  &  & 103 & 0.0131   \\
            4 &57.1  &0  &222.9  &  &  &  &  &  &  &  &  & 419 &  0.009  \\
            5 &61.8  &0  &2.4  &509.7  &  &  &  &  &  &  &  &  637&0.0058    \\
            6 &59.3  &0  &29  &48.8  &947.8  &  &  &  &  &  &  &  1209+& 0.0044   \\
            7 &49.9  &16.2 &51.0  &24.2  &11.1  &1358  &  &  &  &  &  & 741 &0.0039    \\
            8 &39.8  &35.7  &36.3  &9.8  &12.4  &46.2  &1615.4  &  &  &  &  & 959 & 0.0035  \\
            9 &39.2  &33  &40.1  &62.7  &64.1  &13.1  &14.9  & 2326 &  &  &  & 1246+ & 0.0034  \\
            10 &40.9  &21.5  &35.6  &16.4 & 18.9&23.6  &1.1  & 17.5 & 2685.3 &  &  & 1279 &0.0031     \\
            11 &44  &17.5  &38.3  &22.7  &20.4  &28.7  &1.1  &32.3  &44.8  &3487.2  &  &  715&0.0031    \\
            12 &45.5  &5.1  &61.2  &26.9  &2.6  &47.2  &52.5  &81  &75.6  &84.2  &2519.2  & 1370+ & 0.0031  \\

            \end{tabular*}
    {\rule{\temptablewidth}{1pt}}
    \end{center}
\end{table}

% START 2ND SUBSECTION HERE
% increment counter
\addtocounter{Subsection}{1}
% ENTER SUBSECTION TITLE AFTER DOT
\begin{flushleft} \textit{\Alph{Subsection}. A Rational Function Class for the PML Profile}
\end{flushleft}\vspace*{-7pt}

% TYPE THE TEXT OF YOUR 2ND SUBSECTION HERE
First of all, we investigate into (6). Take note that the conditions
$\sigma(0)=0$ and $\sigma'(0)=0$ are satisfied. Let
\begin{equation}
    S_p = ( a_2, a_3, ... , a_p).
\end{equation}
In order to obtain the local optimal values of its coefficients, we
carry out the NM simplex method by using a simple initial value
$S_p=(0,...,0,50)$ for all $p=2,3,...12$, respectively. The result
is summarized in Table 1.

In general, from Table 1, it can be deduced that when $p$ gets
bigger, the number of iteration gets bigger. And $\overline{|R|}$
improves when $p$ increases from 2 to 10, but stopped improving when
$p>10$. The best result comes at $p=10$ and 11 by giving
$\overline{|R|}$=0.0031, only 24\% of the optimal value 0.013
obtained in~\cite{Lu1}. The number of iteration are 1279 and 715,
respectively.

We also observe that: (a) in Table 1, some of the coefficients $a_3,
a_4, ... , a_{p-1}$ are not significant when compared to $a_2$ and
$a_p$; (b) when $\sigma$ gets smaller more rapidly as
$\tau\rightarrow 0$ and bigger more rapidly as $\tau\rightarrow 1$,
the result gets much better.

%define counter for figures
\newcounter{fig}

\addtocounter{fig}{1}
% INCLUDE YOUR PS FIGURE FILE HERE
% Including picture files in Latex is tricky and depends on your graphics package
\vspace*{13pt}\begin{figure}[h]\centering
\includegraphics*[width=0.55\textwidth,angle=0,trim=0 0 10 10, keepaspectratio=true]
    {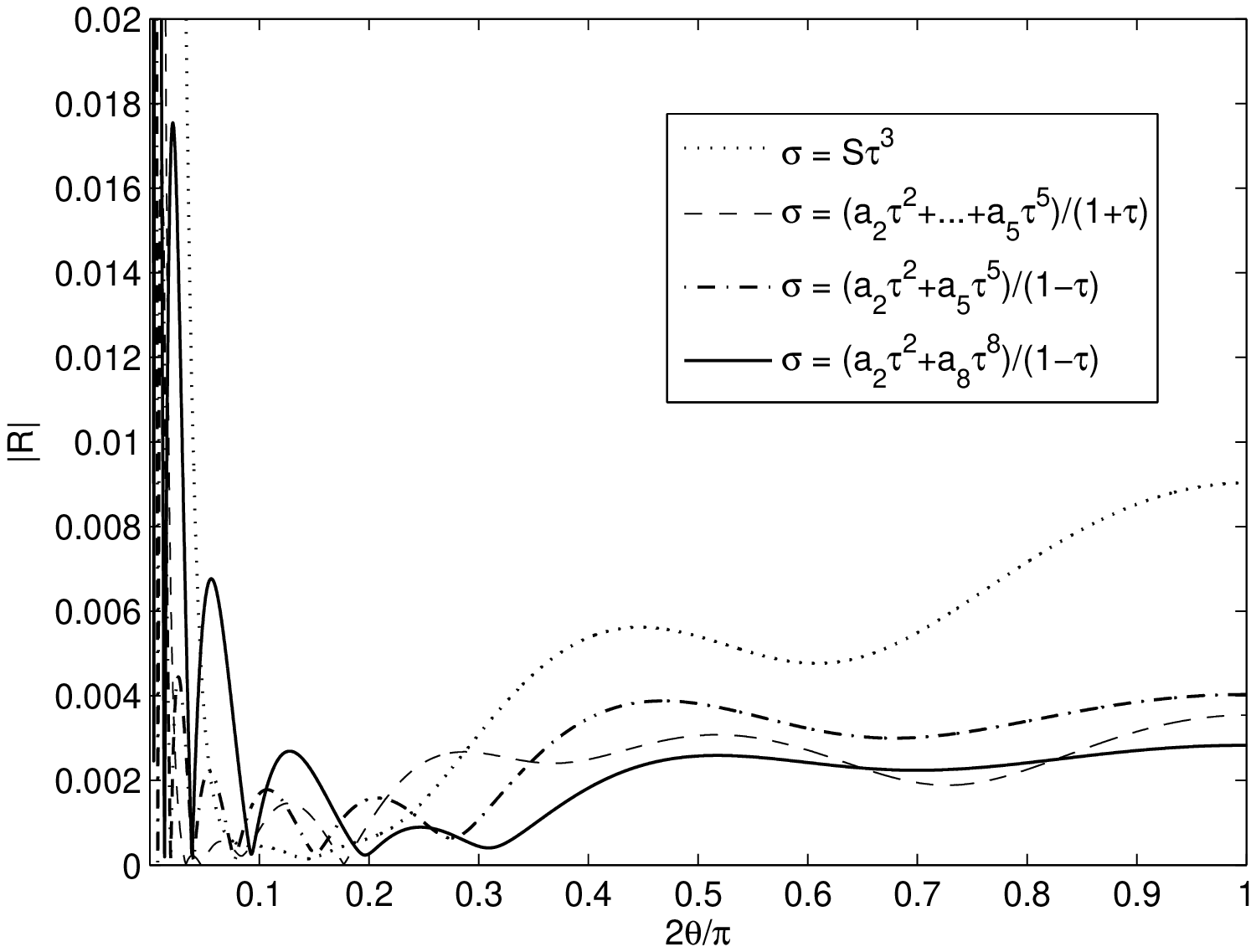}\\
     \centering{Fig. \arabic{fig}. Discrete reflectivity of four different PML profiles using the optimal values listed in Table 1 and Table 2. For the dotted line, $S=100.4$.}
\end{figure}

\addtocounter{fig}{1}
% INCLUDE YOUR PS FIGURE FILE HERE
% Including picture files in Latex is tricky and depends on your graphics package
\vspace*{13pt}\begin{figure}[h]\centering
\includegraphics*[width=0.55\textwidth,angle=0,trim=0 0 10 10, keepaspectratio=true]
    {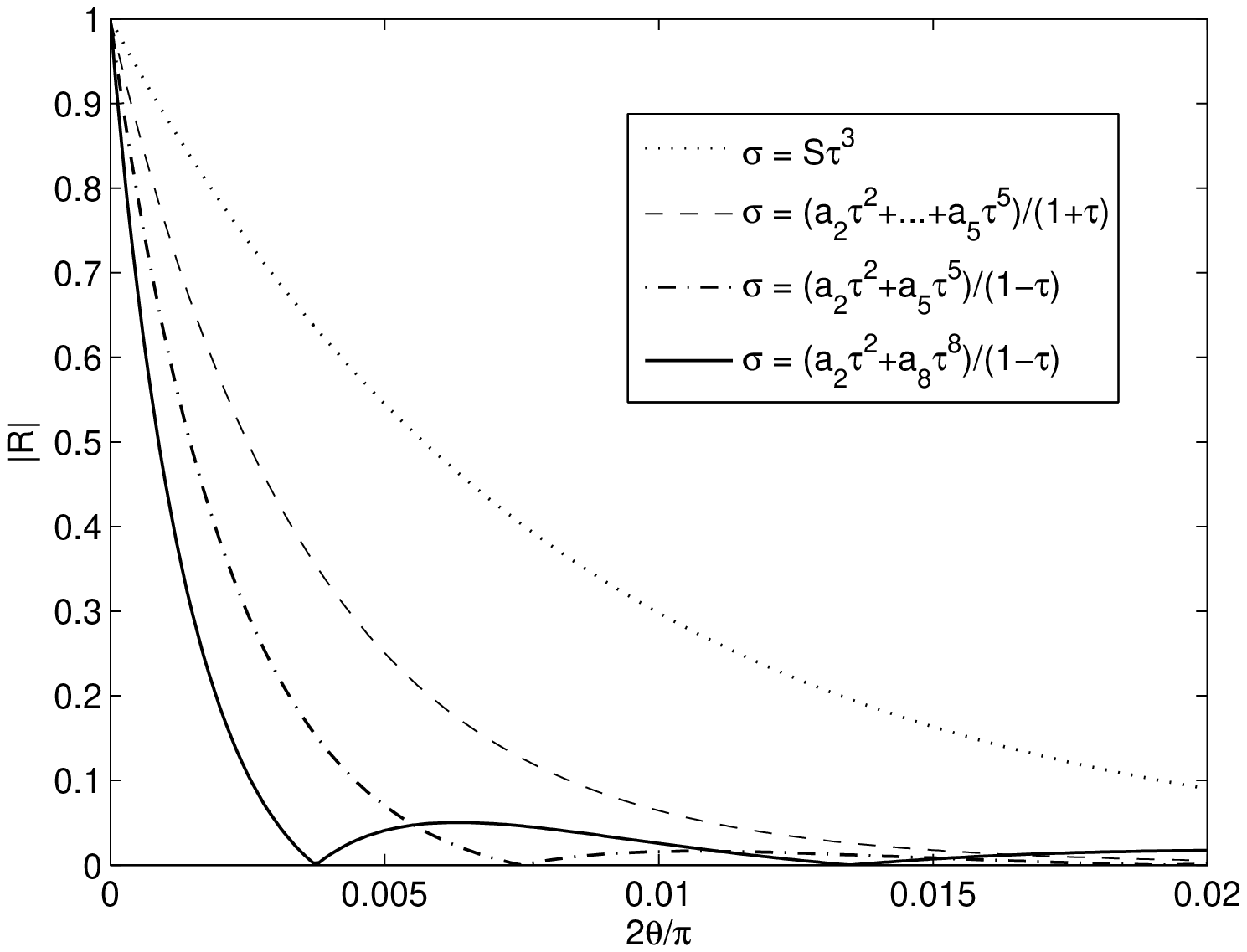}\\
     \centering{Fig. \arabic{fig}. An enlargement of Fig. 1 for small angles.}
\end{figure}

% START 3RD SUBSECTION HERE
% increment counter
\addtocounter{Subsection}{1}
% ENTER SUBSECTION TITLE AFTER DOT
\begin{flushleft} \textit{\Alph{Subsection}. A Simpler Function Class for the PML Profile}
\end{flushleft}\vspace*{-7pt}

Based on the above observations, we try to define a simpler - in
fact, better - function class for $\sigma$: due to (a), let $a_k=0$
for $k=3,...,p-1$ to reduce complexity; due to (b), let the
denominator of $\sigma$ in (6) be $(1-\tau)$, such that
$\sigma\rightarrow 0$ as $\tau\rightarrow 0$ and $\sigma\rightarrow
\infty$ as $\tau\rightarrow 1$. After these changes, the shape of
$\sigma$ satisfies the characteristic in (b) better. Then (6)
becomes
\begin{equation}
    \sigma = \frac{a_2\tau^2+a_p\tau^p}
    {1-\tau}, p\geq2
\end{equation}
where $\sigma(0)=0$ and $\sigma'(0)=0$ still hold.

We carry out the NM simplex method for $\sigma$ defined in (9) again
with the initial value $S_p=(0,...,0,50)$ for all $p=2,3,...12$,
respectively. The local optimal values of its coefficients is listed
in Table 2.

\begin{table}[!h]
\tabcolsep 0pt \caption{Local optimal values of the coefficients of
$\sigma$ defined by (9). Results are derived by the Nelder-Mead
simplex method.} \vspace*{-12pt}
    \begin{center}
    \def\temptablewidth{0.5\textwidth}
    {\rule{\temptablewidth}{1pt}}
        \begin{tabular*}{\temptablewidth}{@{\extracolsep{\fill}}cccccccccccccc}
            $p$ &$a_2$ &$a_p$ &Iterations &$\overline{|R|}$ \\   \hline
            2 &24.9  &      & 21  & 0.0057 \\
            3 &0.0019&28.4  & 44  & 0.0084 \\
            4 &22.5  &14.5  & 109 & 0.0053 \\
            5 &23.6  &35.9  & 91  & 0.0047 \\
            6 &24.4  &76.2  & 136 & 0.0042 \\
            7 &24.3  &113   & 157 & 0.0038 \\
            8 &23.3  &121.3 & 150 & 0.0037 \\
            9 &23.5  &195   & 116 & 0.0037 \\
            10 &23.2 &180.1 & 101 & 0.0038 \\
            11 &23.5 &221.6 & 124 & 0.0039 \\
            12 &23.5 &223.4 & 133 & 0.0041 \\
            \end{tabular*}
    {\rule{\temptablewidth}{1pt}}
    \end{center}
\end{table}

From Table 2, we can see that after simplification, the maximal
number of iteration is reduced to 157. (In Contrast to Table 1, the
maximal number of iteration is more than 1370.) $\overline{|R|}$
gets smaller when $p$ increases from 2 to 8, but bigger when $p>9$.
The best result $\overline{|R|}$=0.0037 occurs at $p=8$ and 9. It is
bigger than $\overline{|R|}$=0.0031 obtained by (6), however,
compared with its fast convergence, (9) is superior in efficiency
and more practical than (6). In practice, we could choose a profile
which has a lower order: for example, $p=5$, the corresponding
$\overline{|R|}$=0.0047.

% START 4TH SUBSECTION HERE
% increment counter
\addtocounter{Subsection}{1}
% ENTER SUBSECTION TITLE AFTER DOT
\begin{flushleft} \textit{\Alph{Subsection}. Comparison}
\end{flushleft}\vspace*{-7pt}

For more details, we plot $\overline{|R|}$ as functions of $\theta$
for four different $\sigma$ in Fig. 1 and Fig. 2. In Fig. 1, we can
see that $\sigma=(a_2\tau^2+a_8\tau^8)/(1-\tau)$ gives the lowest
average of $\overline{|R|}$ for $\theta>0.3\pi/2$. In Fig. 2, for
very small angles ($\theta <0.005\pi/2$), it again distinguishes
itself from the other three by giving the lowest reflectivity.

To illustrate that the NM method gives local minimizers, we plot
$\overline{|R|}$ as multivariable functions of $a_2$ and $a_8$ using
$\sigma=(a_2\tau^2+a_8\tau^8)/(1-\tau)$ in Fig. 3. It can be
observed that at $a_2=23.3$ and $a_8=121.3$, $\overline{|R|}$ almost
reaches its lowest value.

% INITIALIZE SUBSECTION COUNTER
\setcounter{Subsection}{0}

% START A NEW SECTION HERE
% increment counter
\addtocounter{Section}{1}
% ENTER SECTION TITLE AFTER DOT
\begin{center} \textbf{\arabic{Section}. Conclusion}
\end{center}

% TYPE THE TEXT OF YOUR 4RD SECTION HERE
For frequency domain simulations, a procedure for selecting the
optimal PML profile is established based on minimizing the average
discrete reflectivity. By extending the profile to a rational
function class and adopting the Nelder-Mead simplex method to
calculate the profile's coefficients, we reach a better numerical
result. We also provide a simpler profile, which largely conserves
the good properties of the former rational function.

For further improvements, we may use a better function class for the
PML profile, or improve the convergence property of the optimization
method. In addition, we may study the impact of such proposed
profile in a more practical example, for instance, optical wave
propagating along an optical waveguide and hitting the PMLs in
different angles, or how the guided and evanescent waves behave at a
waveguide discontinuity, etc.

\addtocounter{fig}{1}
% INCLUDE YOUR PS FIGURE FILE HERE
% Including picture files in Latex is tricky and depends on your graphics package
\vspace*{13pt}\begin{figure}[h]\centering
\includegraphics*[width=0.5\textwidth,angle=0,trim=0 0 10 10, keepaspectratio=true]
    {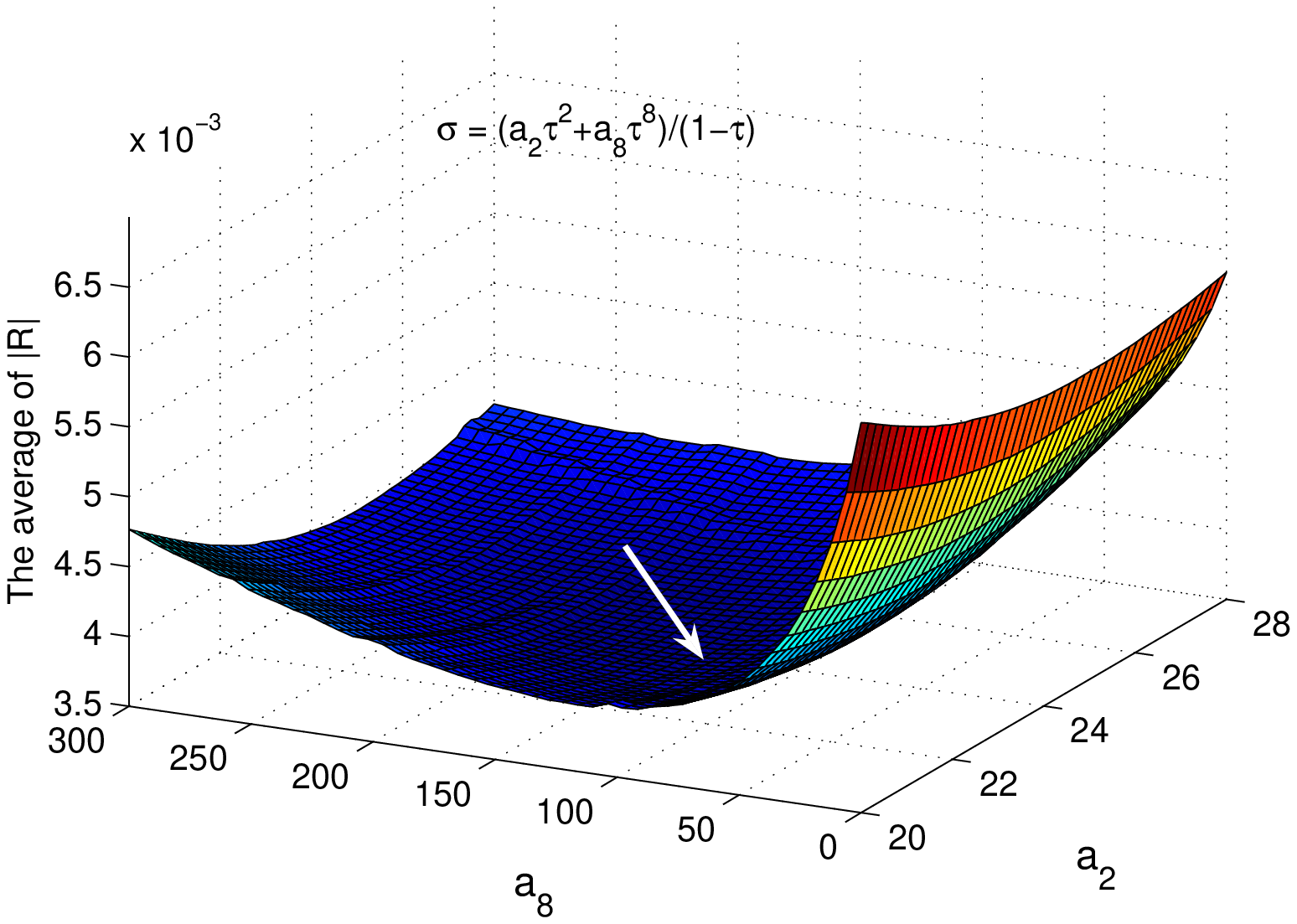}\\
     \centering{Fig. \arabic{fig}. The average of reflectivity as functions of $a_2$ and $a_8$. The white arrow points at $a_2=23.3$ and $a_8=121.3$.}
\end{figure}

\bibliographystyle{IEEEtran}
\bibliography{PMLbib}

\begin{thebibliography}{1}
\providecommand{\url}[1]{#1}
\csname url@rmstyle\endcsname
\providecommand{\newblock}{\relax}
\providecommand{\bibinfo}[2]{#2}
\providecommand\BIBentrySTDinterwordspacing{\spaceskip=0pt\relax}
\providecommand\BIBentryALTinterwordstretchfactor{4}
\providecommand\BIBentryALTinterwordspacing{\spaceskip=\fontdimen2\font plus
\BIBentryALTinterwordstretchfactor\fontdimen3\font minus
  \fontdimen4\font\relax}
\providecommand\BIBforeignlanguage[2]{{%
\expandafter\ifx\csname l@#1\endcsname\relax
\typeout{** WARNING: IEEEtran.bst: No hyphenation pattern has been}%
\typeout{** loaded for the language `#1'. Using the pattern for}%
\typeout{** the default language instead.}%
\else
\language=\csname l@#1\endcsname
\fi
#2}}

\bibitem{Berenger1}
J.~P. Berenger, ``A perfectly matched layer for the absorption of
  electromagnetic waves,'' \emph{Journal of Computational Physics}, vol. 114,
  pp. 185--200, 1994.

\bibitem{Lu1}
Y.~Y. Lu, ``Minimizing the discrete reflectivity of perfectly matched layers,''
  \emph{IEEE Photonics Technology Letters}, vol.~18, no.~3, pp. 487--489, 2006.

\bibitem{YYS1}
D.~Yevick, J.~Yu, and F.~Schmidt, ``Analytic studies of absorbing and
  impedance-mathced boundary layers,'' \emph{IEEE Photonics Technology
  Letters}, vol.~9, no.~1, pp. 73--75, 1997.

\bibitem{LZ1}
Y.~Y. Lu and J.~Zhu, ``Propagating modes in optical waveguides terminated by
  perfectly matched layers,'' \emph{IEEE Photonics Technology Letters},
  vol.~17, no.~12, pp. 2601--2603, 2005.

\bibitem{LRWW1}
J.~C. Lagarias, J.~A. Reeds, M.~H. Wright, and P.~E. Wright, ``Convergence
  property of the nelder-mead simplex method in low dimensions,'' \emph{SIAM
  Journal of Optimization}, vol.~9, no.~1, pp. 112--147, 1998.

\end{thebibliography}

\end{document}